\documentclass[11pt]{amsart}
\usepackage{epsfig}
\usepackage{graphicx, subfigure}
\usepackage{amsmath, amssymb,mathabx,mathrsfs}
\usepackage[toc,page]{appendix}
\usepackage{comment}
\usepackage{hyperref}
\usepackage{color,soul}
\newtheorem{thm}{Theorem}[section]

\theoremstyle{definition}

\theoremstyle{remark}
\newtheorem{rem}[thm]{Remark}
\numberwithin{equation}{section}
\begin{document}
\title[ ]{Regularization of the Hill four-body problem with oblate bodies}

\author{Edward Belbruno}
\address{Yeshiva University, Department of Mathematical Sciences, New York, NY 10016, USA }
\email{edward.belbruno@yu.edu}

\author{Marian Gidea}
\address{Yeshiva University, Department of Mathematical Sciences, New York, NY 10016, USA }
\email{Marian.Gidea@yu.edu}

\thanks{E.B. and M.G. were partially supported by NSF grant  DMS-1814543.}

\author{Wai-Ting Lam}
\address{ Florida Atlantic University, Department of Mathematical Sciences, Boca Raton, FL 33431, USA }
\email{WaiTing.Lam@yu.edu}

\thanks{W-T.L. was partially supported by NSF grant DMS-1814543 and DMS-2138090.}

\begin{abstract}
We consider the  Hill four-body problem where  three oblate, massive bodies form  a  relative equilibrium triangular configuration, and the fourth, infinitesimal body orbits in a neighborhood of the smallest of the three massive bodies.  We regularize collisions between the infinitesimal body and the smallest massive body, via McGehee coordinate transformation. We describe the corresponding collision manifold and show that it undergoes a bifurcation when the oblateness coefficient of the small massive body passes through the zero value.
\end{abstract}

\maketitle
\section{Introduction}

We consider the Hill approximation of the circular restricted four-body problem with oblate bodies, on the motion of an infinitesimal body under the gravitational influence of  three massive  bodies of oblate shapes; the three bodies are assumed to be in a relative equilibrium triangular configuration, and the motion of the infinitesimal body is assumed to take place in a neighborhood of the smallest of the three bodies, which we think of as an asteroid.  See \cite{burgos2020hill}.
The resulting gravitational field in the Hill approximation contains a  non-Newtonian term which depends on the oblateness coefficient of the asteroid.
We use McGehee coordinates to regularize collisions between the infinitesimal body and the asteroid, which amounts to blowing up the collision set to a manifold that captures the dynamics in the singular limit.  (Note that, due to the non-Newtonian term in the potential, the  Levi-Civita regularization does not apply to this setting.)  We describe the collision manifold and the regularized dynamics in a neighborhood of it.
We show that each collision solution is branch regularizable,  and each extension  of a collision solution is a reflection.
We also show that the collision manifold is not block regularizable.
Moreover, we show that the collision manifold undergoes a double saddle-node bifurcation as the oblateness coefficient of the asteroid passes through the zero value. When the shape of the asteroid becomes prolate, no collisions between the infinitesimal body and the asteroid are possible.

The four body system that we consider here can be viewed as a model for the  Sun-Jupiter-Hektor-Skamandrios system; Hektor is a {Jupiter's trojan asteroid}, and Skamandrios is  a moonlet of Hektor. Hektor's shape can be approximated by a dumb-bell figure and has one of the largest oblateness coefficients amongst objects of similar size in the solar system \cite{DESCAMPS2015}.
The moonlet Skamandrios appears to have a complicated orbit, which is close to 1:10 and 2:21 orbit/spin resonances; a
small change could potentially eject the moonlet or make it collide with the asteroid \cite{Marchis}. This justifies our interest in understanding collision orbits.

McGehee coordinate transformation was introduced in \cite{mcgehee1981double} to regularize collisions in  a central force field of the form $U({\bf x})=|{\bf x}|^{-\alpha}$, where ${\bf x}\in\mathbb{R}^2$ and $\alpha>0$. He  also introduced  the concept of branch regularization. A solution is branch regularizable if it has a unique real analytic extension past the collision. Branch regularization concerns the extension of individual solutions.
The concept of block regularization considers  collective extensions of solutions; it was introduced by Easton in \cite{easton1971regularization} who  referred  to it as `regularization by surgery'. A flow is called block regularizable if it is  diffeomorphic to the trivial parallel flow in a deleted neighborhood of the collision set.

{McGehee transformation has been applied to show the existence of ejection-collision orbits, which start and end at a collision.  Llibre showed analytically the existence of ejection-collision orbits in the restricted three-body problem} \cite{llibre1982restricted},
{Lacomba and Llibre showed numerically the existence of transverse
ejection-collision orbits in the Hill problem for some value of the energy} \cite{lacomba1988transversal}, {while} {Delgado-Fern\'andez showed analytically the existence of such orbits  for all sufficiently small energies in} \cite{fernandez1988transversal}. {Other related works include} \cite{devaney1981singularities,pinyol1995ejection,olle2018ejection,alvarez2021ejection}.

McGehee regularization can also be applied to quasi-homogeneous central force fields of the form $U({\bf x})=\gamma_1 |{\bf x}|^{-\alpha_1}+\gamma_1|{\bf x}|^{-\alpha_2}$, with $\gamma_1,\gamma_2,\alpha_1, \alpha_2>0$; see \cite{stoica2000branch}.  Belbruno used McGehee transformation  to regularize collisions with a black hole in \cite{belbruno2011dynamical} in order to establish the relationship between the null geodesic structure of
the Schwarzschild black hole solution, and the corresponding inverse-cubic Newtonian
central force problem. Belbruno and collaborators also used  the McGehee transformation to study the regularizability of the big bang singularity, including the case when random
perturbations modeled by Brownian motion are present in the system \cite{belbruno2013regularizability,xue2014regularization,belbruno2018regularization}.
Other applications of related interest include \cite{diacu2000phase,galindo2014mcgehee,elbialy2009collective,olle2022study}.

A contribution of our work is that we perform McGehee regularization of collisions in a four-body problem (rather than in a  central force field), where the non-Newtonian part of the  gravitational potential is owed to the shape of the body. As a matter of fact, our work assumes a more general setting, of a Hill four-body problem with a general quasi-homogenous potential, which includes the oblateness effect as a particular case.   Another contribution is that we perform a  bifurcation analysis as the oblateness coefficient varies, with the surprising  conclusion that collisions cease to occur as we switch from oblate to prolate shape.

\section{Setup  and main result}
\subsection{Hill four-body problem with oblate bodies}
In this section we describe the Hill approximation of the circular restricted four-body problem with oblate massive bodies. This problem concerns the dynamics of an infinitesimal body (particle) moving in a plane under the gravitational influence of three oblate bodies  of masses $m_1>m_2>m_3$, but without influencing their motion. We refer to these  three bodies as primary, secondary, and tertiary, respectively.
We express the gravitational potential of each body in terms of spherical harmonics  truncated up to second order zonal harmonic, that is,
\begin{equation}\label{eqn:potential}V_i(x_1,x_2,x_3)=\frac{m_i}{r}+\frac{m_i}{r} \left(\frac{R_i}{r}\right)^2 \left(\frac{C^i_{20} }{2}\right)
\left (3 \left (\frac{x_3}{r}\right)^2 -1\right)
\end{equation}
where $r=(x_1^2+x_2^2+x_3^2)^{1/2}$, $R_i$ is the average radius of the $i$-th body,  and the gravitational constant is  normalized  to $1$. The dimensionless quantity $C^i_{20} $ is the coefficient of the zonal harmonic of order $2$, with $C^i_{20} <0$ for an oblate body, $C^i_{20} =0$ for a spherical body,  and
$C^i_{20}>0$ for a prolate body. Further, we denote $C_i=C^i_{20} R_i^2/2$.

\begin{figure}\label{fig:relative_equilibrium}
\includegraphics[width=0.5\textwidth]{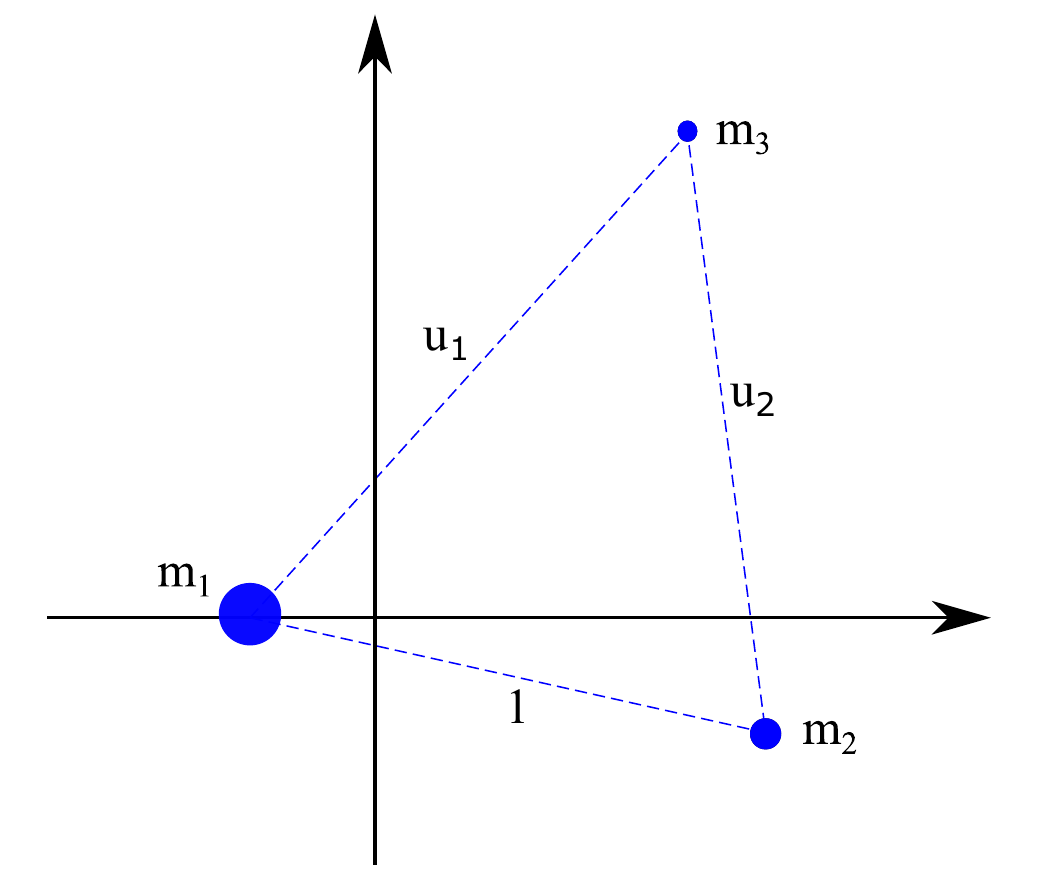}
\caption{Scalene triangle relative equilibrium.}
\end{figure}

For the circular restricted four-body problem, the assumption is that the three massive bodies are in a relative equilibrium configuration, that is, they move on circular orbits around their center of mass while preserving their mutual distances constant over time. In the case when the bodies have no oblateness, the only non-collinear relative equilibrium configuration is the Lagrangian
equilateral triangle. When the bodies are oblate,   the  gravitational
potential is no longer Newtonian, and the relative equilibrium is no longer an equilateral triangle.
It has been shown in \cite{burgos2020hill} that  there is a unique relative equilibrium which is a scalene triangle. Such triangle has the property that the body with the larger $C_i$ is opposite to the longer side of the triangle. We normalize the units of distance so that the distance between $m_1$ and $m_2$ is set to $1$, and we let $u_1$ be the distance from $m_1$ to $m_3$, and $u_2$ be the distance from $m_2$ to $m_3$. See Fig.~\ref{fig:relative_equilibrium}.
The sides $u_1$ and $u_2$ are uniquely determined by the implicit equations
\begin{equation}\label{eqn:u1u2}
 1-3C_{12}=\frac{1}{u_1^3}-\frac{3C_{13}}{u_1^5}  = \frac{1}{u_2^3}-\frac{3C_{23}}{u_2^5},
\end{equation}
where we denote $C_{ij}=C_i+C_j$.

Given such a relative equilibrium configuration,  the motion of the particle in a vicinity of the tertiary is described by the Hamiltonian of the circular restricted four-body problem (see, e.g., \cite{burgos2020hill}). However, the corresponding Hamiltonian equations are difficult to treat analytically. Therefore we consider below the Hill approximation of the circular  restricted four-body
problem. This is derived by  rescaling  the distances  by a factor
of $m_3^{1/3}$,  writing  the associated  Hamiltonian in the
rescaled coordinates as a power series in $m_3^{1/3}$, and
neglecting all the terms of order $O(m_3^{1/3})$  in  the expansion.
The oblateness coefficient $C_i$ also gets rescaled to $c_i=m_3^{-2/3}C^i_{20} R_i^2/2$.
This procedure yields an approximation of the motion of
the particle in an $O(m_3^{1/3})$-neighborhood  of the tertiary,
while the primary and the secondary are `sent to infinity'.
We obtain a much simpler Hamiltonian
{than} the one for the circular restricted four-body problem,
for which the contribution of the primary and of the secondary
to the gravitational potential is given by a quadratic polynomial.
Specifically, the Hamiltonian  of the Hill four-body problem  relative to some convenient co-rotating frame is given by
\begin{equation}\label{eqn:hill_rotated}\begin{split}
H=&\frac{1}{2}(y_1^2+y_2^2+y_3^2)+x_2 y_1-x_1y_2\\
&\quad +\left(\frac{1-\lambda_2}{2}\right)x_1^2+\left(\frac{1-\lambda_1}{2}\right)x_2^2+
\frac{1}{2}\left(\frac{(1-\mu)}{u_1^3}+\frac{\mu}{u_2^3}\right)x_3^2\\
&\quad - \left(\frac{(1-\mu)c_1}{u_1^3}\right)\left(3\left(\frac{x_3}{u_1}\right)^2-1\right)-
\left(\frac{\mu c_2}{u_2^3}\right)\left(3\left(\frac{x_3}{u_2}\right)^2-1\right) \\
&\quad  -\frac{1}{(x_1^2+x_2^2+x_3^2)^{1/2}}- \frac{c_3}{(x_1^2+x_2^2+x_3^2)^{3/2}} \left(\frac{3x_3^2}{ x_1^2+x_2^2+x_3^2 }-
1\right).
\end{split}\end{equation}
where $\mu=\frac{m_2}{m_1+m_2}$, and $\lambda_1$  and $\lambda_2$ are given by the following formulas
\begin{equation}\label{eqn:lambda}\begin{split}
\lambda_1=&\frac{1}{2}\left(2-\frac{2(1-\mu)}{u_1^5}-\frac{2\mu}{u_2^5}+\frac{3(1-\mu)}{u_1^3}
+\frac{3\mu}{u_2^3} -\frac{3}{u_1^3u_2^3}\sqrt{\Delta}\right),\\
\lambda_2=&\frac{1}{2}\left(2-\frac{2(1-\mu)}{u_1^5}-\frac{2\mu}{u_2^5}+\frac{3(1-\mu)}{u_1^3}
+\frac{3\mu}{u_2^3}+\frac{3}{u_1^3u_2^3}\sqrt{\Delta}\right),
\end{split}
\end{equation}
where
\begin{equation*}
\Delta=(\mu u_1^3+(1-\mu)u_2^3)^2-\mu(1-\mu)u_1u_2\left(-u_1^4-u_2^4+2u_1^2+2u_2^2+2u_1^2u_2^2-1\right).
\end{equation*}

When we restrict to the planar problem ($x_3=0$) the Hamiltonian becomes
\begin{equation}\label{eqn:hill_rotated_planar}\begin{split}
H=&\frac{1}{2}(y_1^2+y_2^2)+x_2 y_1-x_1y_2+\left(\frac{1-\lambda_2}{2}\right)x_1^2+\left(\frac{1-\lambda_1}{2}\right)x_2^2\\
&\quad  -\frac{1}{(x_1^2+x_2^2)^{1/2}}+ \frac{c_3}{(x_1^2+x_2^2)^{3/2}},
\end{split}\end{equation}
where the constant terms $\frac{(1-\mu)c_1}{u_1^3}$ and
$\frac{\mu c_2}{u_2^3}$ were dropped, as they do not appear in the Hamiltonian equations.
We note that in the planar problem the oblateness of the primary and the secondary plays no role.

We denote by $\mathbf{M}_h$ the $3$-dimensional energy manifold
\begin{equation}\label{eqn:M}
  \mathbf{M}_h=\{H=h\}.
\end{equation}

\begin{rem}\label{rem:Hektor}
An example of a system that can be modeled by the Hill four-body problem is the Sun-Jupiter-Hektor-Skamandrios system \cite{burgos2020hill}.   Hektor is a Jupiter Trojan, which is approximately located at
Lagrangian point  $L_4$ of the Sun-Jupiter system, thus forming an approximate triangular relative equilibrium configuration with Sun and Jupiter. Hektor is the biggest Jupiter Trojan and
has one of the largest values of the oblateness coefficients among the objects of its size in
in the Solar system.
Hektor's moonlet, Skamandrios, can be viewed as the fourth, infinitesimal body.
In this case, the constants that appear in \eqref{eqn:hill_rotated} are  $c_3=-1.327161\times 10^{-7}$, $\mu=0.0009533386$, $u_1=1- 5.94154\times 10^{-11}$,
$u_2= 1-1.99318\times 10^{-12}$, $\lambda_1= 0.002144$, and $\lambda_2=  2.997855$.
\end{rem}

\subsection{Main result}
The main result of the paper is stated below, and the proof is given in Sections \ref{sec:mcgehee} and \ref{sec:collision_manifold}.
\begin{thm}\label{thm:main}
For the system \eqref{eqn:hill_rotated_planar} with oblate tertiary, i.e., $c_3<0$, each collision solution is branch regularizable,  and each extension  of a collision solution is a reflection. The collision manifold is not block regularizable.

At $c_3=0$ the reduced system of equations associated to the collision manifold undergoes a double saddle-node bifurcation.
For $c_3=0$, the collision manifold is branch and block regularizable.

For the system   \eqref{eqn:hill_rotated_planar} with a prolate  tertiary, i.e., $c_3>0$, there are no collisions.

\end{thm}

The  collision manifold and the corresponding reduced system of equations are described  in Section \ref{sec:collision_manifold}.
%The proof of this result is provided in the following sections.

\section{Branch and block regularization}

We give a brief review of branch and block regularization following \cite{mcgehee1981double}.

For a  differential equation
\begin{equation}\label{eqn:ODE}\dot{\mathbf{x}} = F(\mathbf{x})\end{equation}
with $F$  a real analytic vector field on some open set $U\subseteq \mathbb{R}^n$, and $\dot{}=\frac{d}{dt}$. The standard existence and uniqueness theorem for ODE's gives for each initial condition $\mathbf{x}(0)\in U$ a unique, real analytic solution $\mathbf{x}(t)$ defined on a maximal interval $(t^-,t^+)$ with  $-\infty\leq t^- < 0 < t^+\leq +\infty$.
Solutions for which $-\infty< t^-$ or $t^+<+\infty$ are said to have a singularity  at $t^*=t^-$ or $t^*=t^+$, respectively.

We briefly describe the concept of branch regularization.

If $\mathbf{x}_1(t)$ and $\mathbf{x}_2(t )$ are solutions of \eqref{eqn:ODE}, with $\mathbf{x}_1$ ending
in a singularity at time $t^*$ and $\mathbf{x}_2$ beginning  in a singularity at $t^*$, and there exists  a multivalued analytic complex  function having a branch at $t^*$ and
extending both $\mathbf{x}_1$ and $\mathbf{x}_2$ when we regard the time $t$ as complex,  then  $\mathbf{x}_1$, $\mathbf{x}_2$ are  said to be branch extensions of one another at $t^*$.

\begin{figure}
$\begin{array}{cccc}
\includegraphics[width=0.23\textwidth]{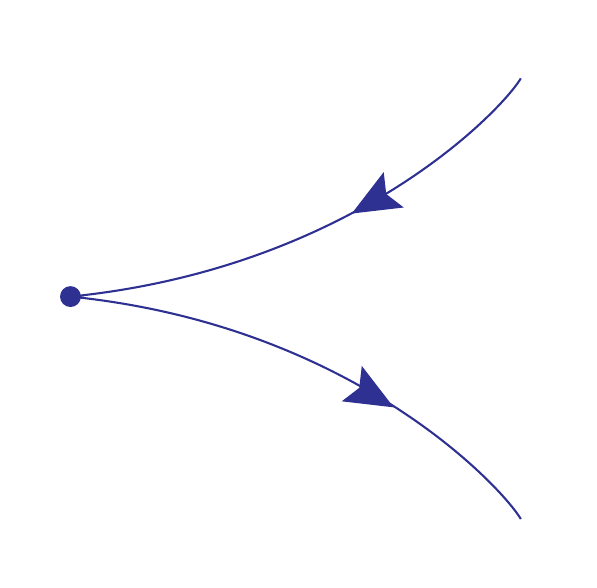}&
\includegraphics[width=0.23\textwidth]{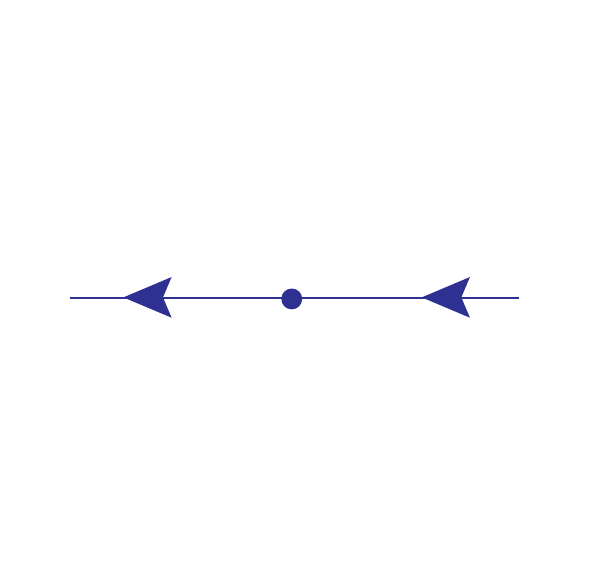}&
\includegraphics[width=0.23\textwidth]{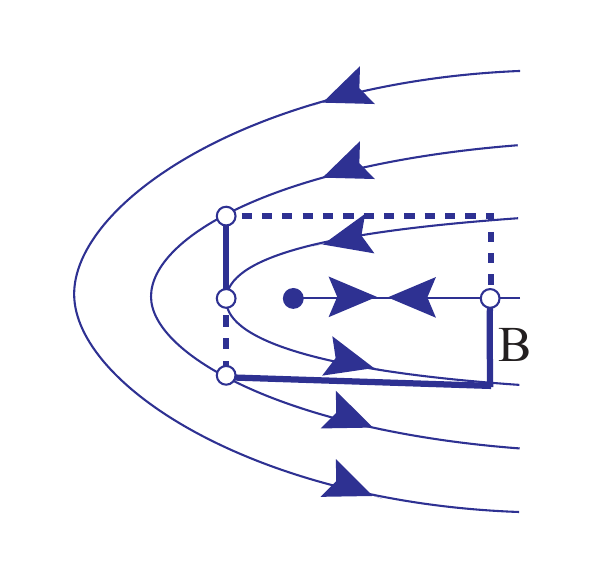}&
\includegraphics[width=0.23\textwidth]{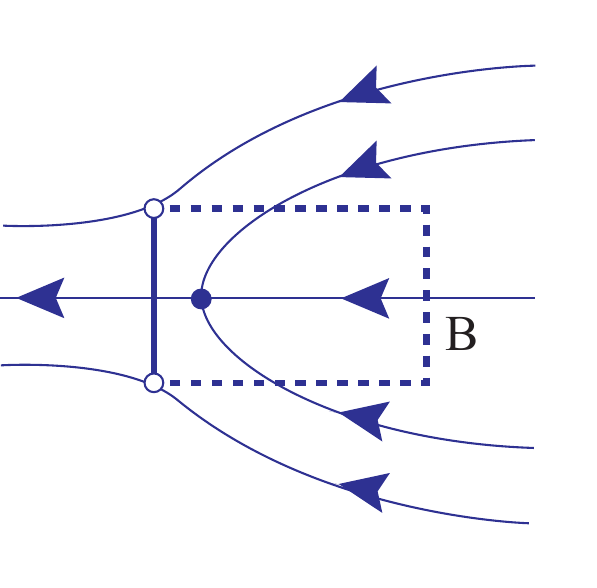}\\
(a) & (b) &(c) &(d)
\end{array}$
\caption{{Different types of singularities: (a) branch regularizable  -- reflection, (b) branch regularizable  -- transmission,  (c) block regularizable -- the sets $a^+$, $a^-$ associated to the  block $B$ are marked by dotted and solid lines, respectively, and the set $t$ is marked by empty circles,  (d) not block regularizable -- same convention as in (c).}}
\label{fig:singularities}\end{figure}

A solution $\mathbf{x}(t)$ of equation\eqref{eqn:ODE} with a  singularity at $t^*$ is said to be branch regularizable at $t^*$ if it has a unique branch
extension at $t^*$.
The extension is called a `reflection'  if  the velocity vector reverses direction at collision, and is called a `transmission'  if  the direction of the velocity vector is preserved at collision. {See Fig.}~\ref{fig:singularities} (a) and (b).

The equation \eqref{eqn:ODE}  is said to be branch regularizable if every solution
is branch regularizable at every singularity.

Now, consider the motion of  a single particle in a potential  field given by
\begin{equation}\label{eqn:potential alpha}U(\mathbf{x})=|\mathbf{x}|^{-\alpha}\textrm{ with $\mathbf{x}\in\mathbb{R}^2$.}\end{equation}
The equation of motion is given by the second order equation \[\ddot{\mathbf{x}}=\nabla U(\mathbf{x}),\] or equivalently, by the first order system
\begin{equation*}
\left\{
\begin{split}
\dot{\mathbf{x}}=&\mathbf{y},\\
\dot{\mathbf{y}}=&-\alpha|\mathbf{x}|^{-\alpha-2} \mathbf{x}.
\end{split}
\right.
\end{equation*}

Let  $\beta=\frac{\alpha}{2}$ and  $\gamma=\frac{1}{1+\beta}=\frac{2}{2+\alpha}$.

We recall the following result from \cite{mcgehee1981double}:
\begin{thm}\label{thm:branch}

A collision solution for the potential \eqref{eqn:potential alpha} is branch regularizable if and only if $\gamma=\frac{p}{q}$ with $p<q$ positive integers, $\textrm{gcd}(p,q)=1$, and  $q$ odd.

Moreover, if $p$ is even the extension solution is a `reflection', and when $p$ is odd the extension solution is a `transmission'.

\end{thm}

In \cite{stoica2000branch} this result has been extended for quasi-homogeneous potentials of the form
\begin{equation}\label{eqn:quasi}
U(\mathbf{x})=\gamma_1|\mathbf{x}|^{-\alpha_1}+\gamma_2|\mathbf{x}|^{-\alpha_2}\textrm { with } \mathbf{x}\in\mathbb{R}^2,
\end{equation} where $\gamma_1,\gamma_2>0$, $\alpha_2>\alpha_1>0$.

\begin{thm}\label{thm:branch_2}
A collision solution for the potential \eqref{eqn:quasi} is branch regularizable if and only if both
\[\frac{2}{2+\max(\alpha_1,\alpha_2)}, \textrm { and } \frac{\min(\alpha_1,\alpha_2)}{2+\max(\alpha_1,\alpha_2)}\]
are of the form $\frac{p}{q}$ with $p<q$ positive integers, $\textrm{gcd}(p,q)=1$, and  $q$ odd.
\end{thm}

We now describe the concept of block regularization. Denote by $\phi^t=\phi(\cdot,t)$ the flow of \eqref{eqn:ODE}.

A compact invariant set  $\mathbf{N}\subseteq U$ is called isolated if there exists an open set $V\subseteq U$ -- referred to as an isolating neighborhood --  such that $\mathbf{N}\subset V$  is the maximal invariant subset of $V$.

Let $\mathbf{B}\subseteq U$ be a compact set with non-empty interior, and assume that the boundary $\mathbf{b}=\partial \mathbf{B}$   of $\mathbf{B}$ is a smooth submanifold. Define
\begin{equation*}
  \begin{split}
     \mathbf{b}^+=&\{ \mathbf{x}\in \mathbf{b} \,|\, \phi( \mathbf{x},(-\varepsilon,0))\cap \mathbf{B}= \emptyset, \textrm{ for some }\varepsilon>0\},\\
     \mathbf{b}^-=&\{ \mathbf{x}\in \mathbf{b} \,|\,  \phi( \mathbf{x},(0,\varepsilon))\cap \mathbf{B}= \emptyset, \textrm{ for some }\varepsilon>0\},\\
     \mathbf{t}=&\{ \mathbf{x}\in \mathbf{b}     \,|\,  \dot{\phi}(\mathbf{x},0)\textrm  { is tangent to } \mathbf{b} \}.
  \end{split}
\end{equation*}

The set $\mathbf{B}$  is called an isolating block if $ \mathbf{b}^+\cap \mathbf{b}^-= \mathbf{t}$.

If $\mathbf{N}$ is an isolated invariant set, we say that $\mathbf{B}$ isolates $\mathbf{N}$ if the interior set $\textrm{Int}(\mathbf{B})$ of $\mathbf{B}$  is an isolating neighborhood for $\mathbf{N}$.
For every isolated invariant set $\mathbf{N}$ there exists an isolating block which isolates $\mathbf{N}$. If $\mathbf{B}$ is an isolating block, then there exists an isolated invariant set $\mathbf{N}$ (possibly empty) which is isolated by $\mathbf{B}$. See \cite{conley1971isolated}.

The asymptotic sets to $\mathbf{N}$ are defined by
\begin{equation*}
  \begin{split}
     \mathbf{a}^+=&\{ \mathbf{x}\in \mathbf{b}^+ \,|\, \phi( \mathbf{x},(0,+\infty))\subset \mathbf{B}\}\\
     \mathbf{a}^-=&\{ \mathbf{x}\in \mathbf{b}^-\,|\,  \phi( \mathbf{x},(-\infty,0))\subset \mathbf{B}\}\\
  \end{split}
\end{equation*}

The map across the block is defined as
\begin{equation*}
  \begin{split}\Phi: \mathbf{b}^+\setminus \mathbf{a}^+ \to \mathbf{b}^-\setminus \mathbf{a}^-,\\
\Phi (\mathbf{x})=\phi(\mathbf{x},T(\mathbf{x})),
  \end{split}
\end{equation*}
where $T(\mathbf{x})=\inf\{t>0\,|\, \phi(\mathbf{x},t)\not\in\mathbf{B}\}$ is the time spent inside the block.

If $\mathbf{B}$ is an isolating block, then the application
$\Phi$  is a diffeomorphism. See \cite{conley1971isolated}.

An isolating block $\mathbf{B}$ is said to be trivializable if the map $\Phi$ extends
uniquely to a diffeomorphism from $ \mathbf{b}^+$ to $\mathbf{b}^-$.

The theory of isolating blocks can be applied to singularities by essentially replacing the role of an isolated  invariant set $\mathbf{N}$ as above with the set of singularities, as we shall see below.

In Section \ref{sec:mcgehee} we will see that, going through regularized coordinates and time rescaling,  the set of singularities for \eqref{eqn:hill_rotated}, which consists of the origin, gets transformed into an invariant set, which is in fact a manifold (referred to as a collision manifold).

Let $F(\mathbf{x})$ be a vector field defined on $U\setminus \mathbf{N}$, where $\mathbf{N}$ is a compact set representing  the singularities of the vector field.  Let $\mathbf{B}\subseteq U$ be compact set with non-empty interior, such that  $\mathbf{b}=\partial \mathbf{B}$   is a smooth submanifold, and with $\mathbf{b}\cap \mathbf{N}=\emptyset$. Define the subsets $\mathbf{b}^+,\mathbf{b}^-\subset \mathbf{b}$ in the same way as above.
Under these conditions, the definition of an isolating block is the same as before.

The orbit through a point $\mathbf{x}$ is defined by
\[O(\mathbf{x}) =\{\phi(\mathbf{x}, t) \,|\, \phi(\mathbf{x}, t)\textrm{ is defined }\}.\]
That is, there are no invariant sets in $\mathbf{B}$.

An isolating block $\mathbf{B}$ is said to isolate the singularity set $\mathbf{N}$ if
$\mathbf{N} \subset\textrm{Int}(\mathbf{B})$  and if $O(\mathbf{x})\not\subset \mathbf{B}$ for all $\mathbf{x}\in \mathbf{B}\setminus\mathbf{N}$.

The asymptotic sets to $\mathbf{N}$ are defined by
\begin{equation*}
  \begin{split}
     \mathbf{a}^+=&\{ \mathbf{x}\in \mathbf{b}^+ \,|\, \phi( \mathbf{x},t)\in \mathbf{B} \textrm { for all $t\geq 0$ for which $ \phi( \mathbf{x},t)$ is defined} \}\\
     \mathbf{a}^-=&\{ \mathbf{x}\in \mathbf{b}^-\,|\,  \phi( \mathbf{x},t)\in \mathbf{B} \textrm { for all $t\leq 0$ for which $ \phi( \mathbf{x},t)$ is defined}\}.
  \end{split}
\end{equation*}
We define the map across the block $\Phi: \mathbf{b}^+\setminus \mathbf{a}^+ \to \mathbf{b}^-\setminus \mathbf{a}^-$ as before.

The singularity set $\mathbf{N}$ is said to be block regularizable if there exists a trivializable block $\mathbf{B}$ which isolates $\mathbf{N}$.
 {See Fig.}~\ref{fig:singularities} (c) and (d).

Regarding block regularization we  recall the following result from \cite{mcgehee1981double}:
\begin{thm}\label{thm:block}
A collision set  for the potential \eqref{eqn:potential alpha} is block regularizable if and only if $\beta=1-\frac{1}{n}$ for $n$ positive integer.
\end{thm}

\section{McGehee transformation}
\label{sec:mcgehee}
We rewrite the Hamiltonian \eqref{eqn:hill_rotated_planar} in a simpler form
\begin{equation}\label{eqn:hill_rotated_simplified}\begin{split}
H=&\frac{1}{2}(y_1^2+y_2^2)+x_2 y_1-x_1y_2+Ax_1^2+B x_2^2\\
&\quad  -\frac{1}{(x_1^2+x_2^2)^{\nu/2}}- \frac{c}{(x_1^2+x_2^2)^{\alpha/2}},
\end{split}\end{equation}
where $\dot{}=\frac{d}{dt}$, and $1\leq \nu< \alpha$.  The corresponding potential is quasi-homogeneous.

In the case of the potential \eqref{eqn:hill_rotated_planar}, we have \begin{equation}\label{eqn:AB}\begin{split}\nu=1,\quad
\alpha=3,\quad
A= \frac{1-\lambda_2}{2},\quad
B=\frac{1-\lambda_1}{2},\quad\lambda_1, \lambda_2>0,\quad
c=-c_3.\end{split}\end{equation}

We identify $x,y\in\mathbb{R}^2$ with the complex numbers $x_1+ix_2$, $y_1+iy_2$, respectively.
The corresponding Hamilton equations are
\begin{equation}\label{eqn:ham_eqn_0}\begin{split}
\dot x =&\frac{\partial H}{\partial y}=y-ix,\\
\dot y =&-\frac{\partial H}{\partial x}= -\frac{\nu x}{|x|^{\nu+2}}-\frac{\alpha c x}{|x|^{\alpha+2}}-iy-Tx,
\end{split}\end{equation}
where  $T$ is the real-linear transformation given by $T(x_1+ix_2)=2Ax_1+2Bx_2 i$, and $|x|=(x_1^2+x_2^2)^{1/2}$.

We perform a coordinate change to new real coordinates $(r,\theta, v,w)$, with $r>0$ and $\theta\in\mathbb{T}^1$, defined as follows
\begin{equation}\label{eqn:mcgehee}\begin{split}
x =& r^\gamma e^{i\theta},\\
y =&   r^{-\gamma\beta}(v+iw)e^{i\theta},
\end{split}\end{equation}
where
\[ \beta=\frac{\alpha}{2}, \textrm{ and } \gamma=\frac{1}{\beta+1}=\frac{2}{\alpha+2}.\]
Writing \eqref{eqn:mcgehee} in terms of components we have
\begin{equation}\label{eqn:mcgehee_comp}\begin{split}
  x_1=&r^\gamma\cos\theta,\\
  x_2=&r^\gamma\sin\theta,\\
  y_1=&r^{-\gamma\beta}(v\cos\theta-w\sin\theta),\\
  y_2=&r^{-\gamma\beta}(v\sin\theta+w\cos\theta).
\end{split}
\end{equation}
The new coordinates $(r,\theta, v,w)$ in terms of the old coordinates $(x_1,x_2,y_1,y_2)$ are given by
\begin{equation}\label{eqn:mcgehee_comp_inv}\begin{split}
  r=&|x|^{\frac{1}{\gamma}},\\
  \theta=&\arg(x),\\
  v=&r^{\gamma\beta}(y_1\cos\theta+y_2\sin\theta),\\
  w=&r^{\gamma\beta}(-y_1\sin\theta+y_2\cos\theta).
\end{split}
\end{equation}

The new coordinates are known as the McGehee coordinates \cite{mcgehee1981double}.
We rewrite the Hamiltonian equations \eqref{eqn:ham_eqn_0} in the new coordinates and equate the real and imaginary parts on the two sides.
From
\begin{equation*}\begin{split}
\dot{x}=&\left[\gamma r^{\gamma-1}\dot{r}+i r^\gamma \dot{\theta}\right]e^{i\theta},\\
  \frac{\partial H}{\partial y}=&\left[r^{-\beta\gamma}(v+iw)-i r^\gamma \right] e^{i\theta},\\
\dot{y}=&\left[-\beta\gamma r^{-\beta\gamma-1}\dot {r}(v+iw)+r^{-\beta\gamma}(\dot {v}+i\dot{w})+r^{-\beta\gamma}(v+iw)i\dot{\theta}\right]e^{i\theta}\\
  =&\left[-\beta r^{-\beta\gamma-1}v(v+iw)+r^{-\beta\gamma}(\dot {v}+i\dot{w})+r^{-\beta\gamma}(v+iw)i(r^{-1}w-1)\right]e^{i\theta},\\
 - \frac{\partial H}{\partial x}=&\left [
-\frac{\nu r^\gamma}{r^{\gamma(\nu+2)}}-\frac{\alpha c r^\gamma}{r^{\gamma(\alpha+2)}}
 - i\frac{v+iw}{r^{\gamma\beta}}\right]e^{i\theta}
 -2Ar^\gamma\cos\theta-i2Br^\gamma\sin\theta\\
 &=\left [-\nu r^{\gamma(-1-\nu)}-\alpha c r^{\gamma(-1-\alpha)}
- i(v+iw)r^{-\gamma\beta}\right.\\
 &\quad\left.+2r^\gamma(-A\cos^2\theta-B\sin^2\theta)
 +i 2 r^\gamma(A-B)\sin\theta\cos\theta  \right]e^{i\theta}
\end{split}\end{equation*}
we obtain
\begin{equation}\label{eqn:ham_eqn_1}\begin{split}
\dot{r}=&(\beta+1)v,\\
\dot{\theta}=&r^{-1}w-1,\\
\dot{v}=&\frac{\beta v^2+w^2-\alpha c}{r}-\frac{\nu}{r^{\gamma(\nu+2)-1}}-2Ar\cos^2\theta-2Br\sin^2\theta,\\
\dot{w}=&\frac{(\beta-1)vw}{r}+2(A-B)r\sin\theta\cos\theta.
\end{split}
\end{equation}
In the above, after equating $\dot x=\frac{\partial H}{\partial y}$ we obtain $\dot{r}= (\beta+1)v$ and $\dot{\theta}=r^{-1}w-1$, which we substitute in the  equation for $\dot y$. We also use  that $\frac{1-\gamma}{\gamma}=\beta$, $-\gamma\beta=\gamma-1$, and $\alpha=2\beta$. The fact that $T$ is real-linear transformation  but not-complex linear is taken into account when factoring out $e^{i\theta}$ in the equation for $ - \frac{\partial H}{\partial x}$  by expressing $Tx=(Tx  e^{-i\theta})e^{i\theta}$.

The equations \eqref{eqn:ham_eqn_1} have a singularity at $r=0$.
We remove the singularity by introducing a new time parameter $\tau$ given by
\begin{equation}\label{eqn:time}
dt=r d\tau.
\end{equation}
The equations  \eqref{eqn:ham_eqn_1}  expressed in terms of the new time $\tau$ become
\begin{equation}\label{eqn:ham_eqn_2}\begin{split}
r'=&(\beta+1)vr,\\
\theta'=&w-r,\\
v'=&( \beta v^2+w^2- \alpha c)-\nu r^{2-\gamma(\nu+2)}-2Ar^2\cos^2\theta-2Br^2\sin^2\theta,\\
w'=&(\beta-1)vw+2(A-B)r^2\sin\theta\cos\theta,
\end{split}
\end{equation}
where ${}'=\frac{d}{d\tau}$. Since $\nu<\alpha$, we have that $2-\gamma(\nu+2)>0$. Thus, the obtained differential equations have no singularity at $r=0$; the singularity has been `removed'. We also note that the terms $\nu r^{2-\gamma(\nu+2)} -2Ar^2\cos^2\theta-2Br^2\sin^2\theta$ and $2(A-B)r^2\sin\theta\cos\theta$ tend to $0$ as $r\to 0$, so they can be neglected for $r$ sufficiently small.

The energy condition $H=h$ in the new coordinates, when we use \eqref{eqn:mcgehee_comp}, becomes
\begin{equation}\label{eqn:energy_h}\begin{split}
&\frac{1}{2}r^{-2\gamma\beta}(v^2+w^2)\\
&+r^{\gamma(1-\beta)} \sin\theta (v\cos\theta-w\sin\theta)-r^{\gamma(1-\beta)}\cos\theta (w\cos\theta+v\sin\theta)\\
&+Ar^{2\gamma}\cos^2\theta+Br^{2\gamma}\sin^2\theta -r^{-\gamma\nu}-cr^{-\gamma\alpha}=h,
\end{split}\end{equation}
which, after multiplying both sides by $r^{2\gamma\beta}=r^{2-2\gamma}$ yields
\begin{equation}\label{eqn:h}
\begin{split}
\frac{v^2+w^2-2c}{2}-rw+r^2(A\cos^2\theta+B\sin^2\theta)-r^{2-\gamma(\nu+2)}=r^{2-2\gamma}\, h.
\end{split}
\end{equation}

{We define the \emph{energy manifold} $\mathbf{M}_h$ as the set of points $(r,\theta,v,w)$ satisfying}~\eqref{eqn:h}.
When $r=0$ the energy condition  \eqref{eqn:h} reduces to
\begin{equation}\label{eqn:h_r_0}
\begin{split}
v^2+w^2-2c=0.
\end{split}
\end{equation}

\begin{rem}
{In the case of the potential} \eqref{eqn:potential alpha}, {one obtains a system of} $4$-{equations  similar to} \eqref{eqn:ham_eqn_2}:
\begin{equation}\label{eqn:mcgehee_eqn}\begin{split}
r'=&(\beta+1)vr,\\
\theta'=&w,\\
v'=&\beta(v^2-2)+w^2,\\
w'=&(\beta-1)vw,
\end{split}
\end{equation}
This   system is partially decoupled  -- the first two equations are determined by the last two
equations. Also, the energy manifold $\mathbf{M}_h$ projects onto $\{v^2+w^2<1\}$ when $h<0$,
onto $\{v^2+w^2>1\}$ when $h>0$, and onto $\{v^2+w^2=1\}$ when $h=0$.
See\cite{mcgehee1981double}.

{In the case of our system} \eqref{eqn:ham_eqn_2} {the second equation  has an extra term owed to the Coriolis effect in} \eqref{eqn:hill_rotated}, {the third equation has extra terms owed to oblateness and to the effects of the primary and secondary, and the fourth equation has an extra term owed to the effect of the primary and secondary.}

Also,   the system \eqref{eqn:ham_eqn_2} is fully coupled, and there is no obvious relation between the regions bounded by $v^2+w^2=2c$ in the $(v,w)$-plane and the energy $h$.

Writing  the energy condition  \eqref{eqn:h} as
\begin{equation*}\begin{split}
& \frac{v^2+w^2-2c}{2}=r^{2-2\gamma}\, h +rw-r^2(A\cos^2\theta+B\sin^2\theta)+r^{2-\gamma(\nu+2)}, \textrm{ or }\\
& h= \frac{1}{r^{2-2\gamma}}\left(\frac{v^2+w^2-2c}{2}\right) - \frac{1}{r^{1-2\gamma}}w +r^{2\gamma}  (A\cos^2\theta+B\sin^2\theta) - \frac{1}{r^{\nu\gamma}},
\end{split}
\end{equation*}
{we see that for $r\ll 1$ the sign of $h$ is the same as the sign of $v^2+w^2-2c$.  Thus, the  points in $\mathbf{M}_h$ with $h>0$ and $r\approx 0$ project onto $\{v^2+w^2>2c\}$ and the  points in $\mathbf{M}_h$ with $h<0$ and $r\approx 0$  project onto $\{v^2+w^2<2c\}$.}

%by neglecting the $O(r^2)$ terms we can obtain that the points in $\mathbf{M}_h$ with $h>0$, $w>0$ and $r$ sufficiently small project onto $\{v^2+w^2>2c\}$.

%Also, provided that $\alpha>2(\nu+1)$, we can neglect $r^{2-\gamma(\nu+2)}$ and obtain that
 %the points in $\mathbf{M}_h$ with $h<0$, $w<0$ and $r$ sufficiently small project onto $\{v^2+w^2<2c\}$.
\end{rem}

\section{ Equilibrium points and Hill regions}

{A straightforward computations shows that} \eqref{eqn:ham_eqn_2}  {has $6$ equilibrium points, which  in terms of the
$(r,\theta,v,w)$-coordinates are given by}
\begin{equation}\label{eqn:equilibria}
\begin{split}
\mathscr{E}_{\pm}=&(0,\theta_0,\pm \sqrt{2c},0),\\
\mathscr{E}_{1}=&(r_{1},0 , 0,r_{1}),\quad \mathscr{E}_{2}=(r_{1},\pi , 0,r_{1}),\\
\mathscr{E}_{3}=&(r_{2},\pi/2 , 0,r_{2}),\quad \mathscr{E}_{4}=(r_{2},3\pi/2 , 0,r_{2}).
\end{split}
\end{equation}
{for arbitrary $\theta_0\in \mathbb{T}^1$,
$r_{1}$ being the solution  of}
\begin{equation}\label{eqn:E12} r^2(1-2A)-\nu r^{2-\gamma(\nu+2)}-\alpha c=0, \end{equation}
{and $r_{2}$ being the  solution  of}
\begin{equation}\label{eqn:E34} r^2(1-2B)-\nu r^{2-\gamma(\nu+2)}-\alpha c=0 .\end{equation}
{Note that} \eqref{eqn:E12}  and \eqref{eqn:E34} {have unique solutions.}
%{In the case when $\nu=1$ and $\alpha=3$, by } \eqref{eqn:AB} {we have $1-2A=\lambda_2>0$ and $1-2B=\lambda_1>0$, and by} \eqref{eqn:mcgehee} {we have $r^\gamma=\|x\|$,  the
{The equilibrium points $\mathscr{E}_1, \mathscr{E}_2, \mathscr{E}_3, \mathscr{E}_4$ are the same as the  $x$- and $y$-equilibrium points for the Hill four-body problem in} \cite{burgos2020hill}, respectively ({referred to as $L_1, L_2,L_3,L_4$ in} \cite{burgos2016families}).
{The points $\mathscr{E}_1, \mathscr{E}_2$ are of center-saddle type; the points $\mathscr{E}_3, \mathscr{E}_4$ are of center-center type  provided that $\mu$ is less than some critical value $\mu_{\textrm{cr}}$.
On the other hand,  $\mathscr{E}_{\pm}$ form circles of equilibrium points.
The eigenvalues at each point of $\mathscr{E}_{\pm}$ are} \[0, \pm(\beta+1)\sqrt{2c},\pm 2\beta\sqrt{2c},\pm(\beta-1)\sqrt{2c}.\]
{The circle $\mathscr{E}_+$ has a $4$-dimensional unstable manifold and the circle $\mathscr{E}_-$ has a $4$-dimensional stable manifold, which necessarily coincide.}
\begin{figure}
$\begin{array}{cc}
\includegraphics[width=0.4\textwidth]{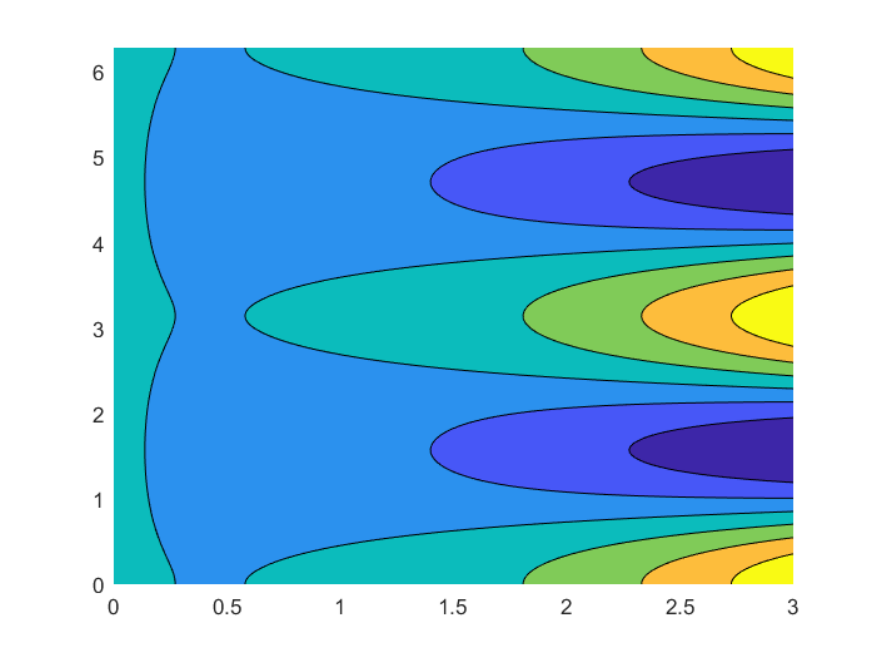}&
\includegraphics[width=0.4\textwidth]{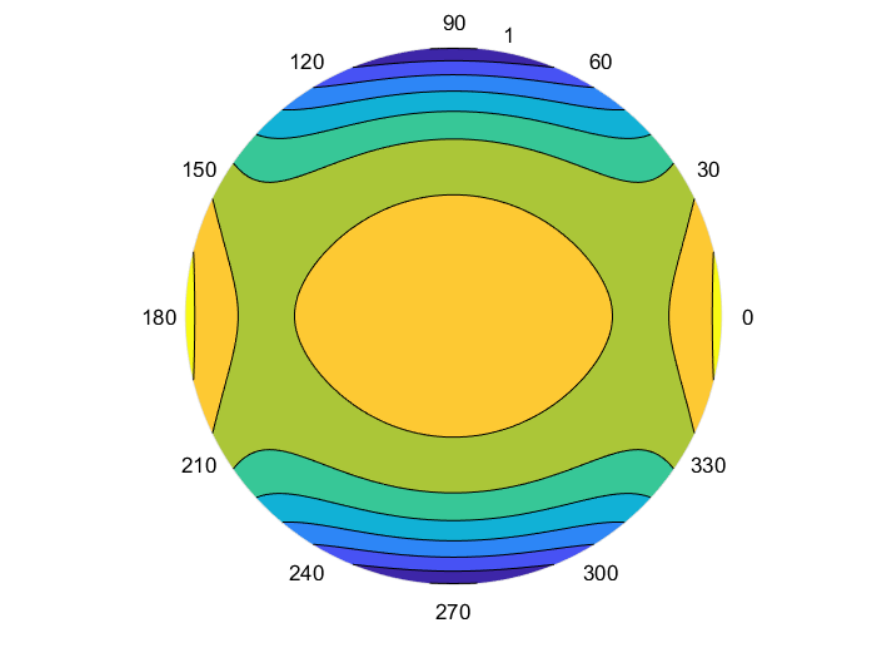}\\
\includegraphics[width=0.4\textwidth]{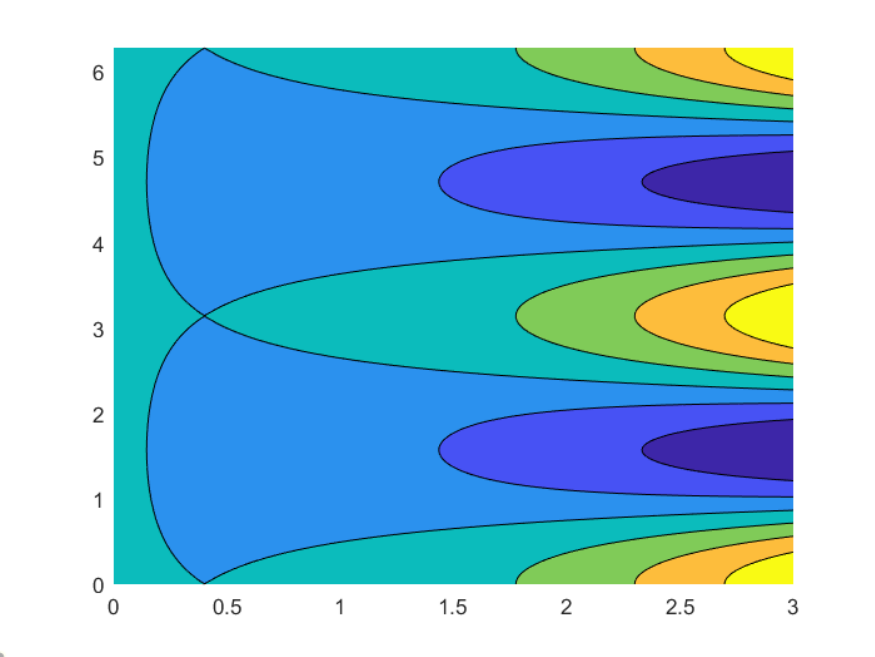}&
\includegraphics[width=0.4\textwidth]{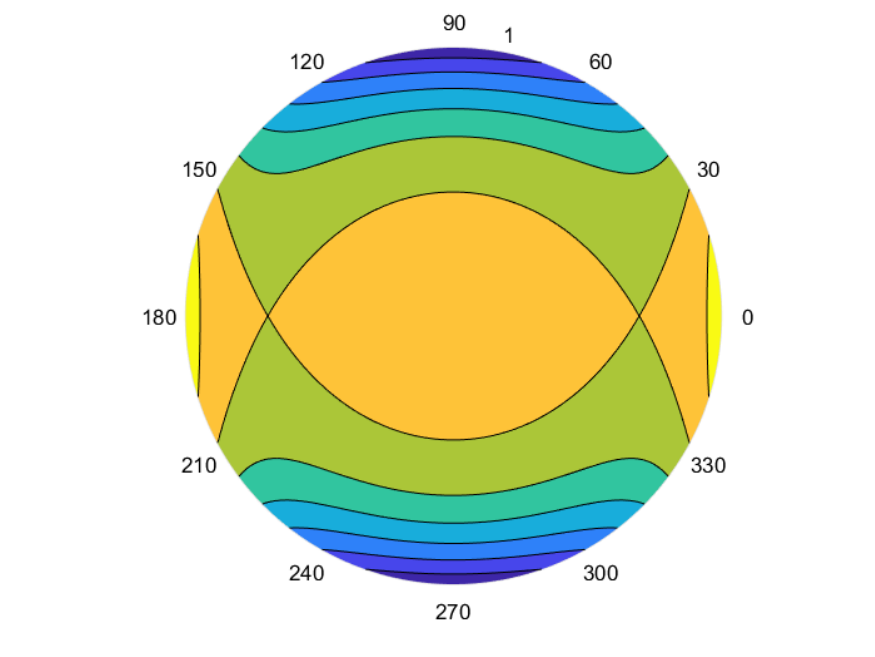}\\
\includegraphics[width=0.4\textwidth]{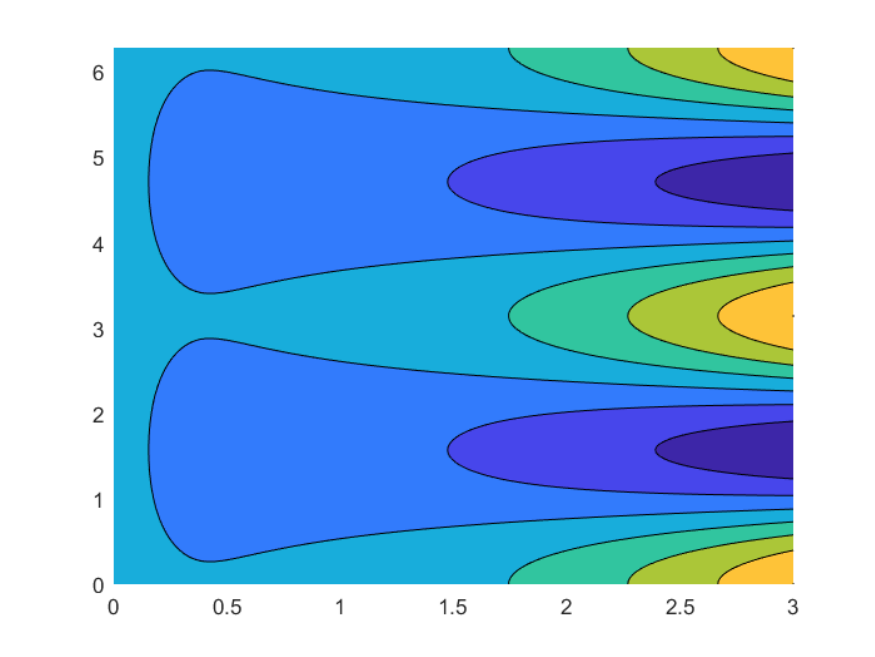}&
\includegraphics[width=0.4\textwidth]{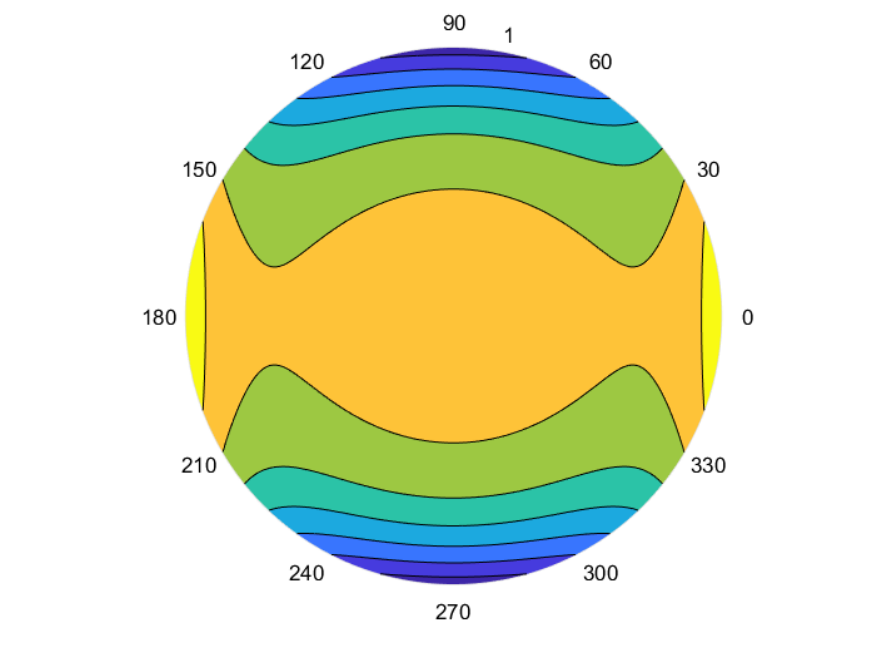}\\
(a) & (b)
\end{array}$
\caption{Hill region for  energy levels below, at, and above that of $\mathscr{E}_1, \mathscr{E}_2$ in McGehee coordinates:  (a) in $(r,\theta)$-coordinates, (b) in $(x,y)$-coordinates. }
\label{fig:Hill_region_0}\end{figure}

{The effective potential for the system} \eqref{eqn:hill_rotated_simplified} {is}
\[\Omega(x_1,x_2)=\frac{1}{2}\left(\lambda_2x_1^2+\lambda_1x_2^2\right)+\frac{1}{(x_1^2+x_2^2)^{\nu/2}}+\frac{c}{(x_1^2+x_2^2)^{\alpha/2}},
\]
{which, written in McGehee coordinates becomes}
\[\Omega(r,\theta)=\frac{1}{2}\left(\lambda_2r^{2\gamma}\cos^2(\theta)+\lambda_1r^{2\gamma} \sin^2(\theta)\right)+
r^{-\nu\gamma}+{c}r^{-\alpha\gamma}.
\]
{Then the  Hill region for an energy level $h$, defined as the projection on the energy manifold onto configuration space, represents  the region of possible motions,  and is given by}
\[\{(r,\theta)\,|\, \Omega(r,\theta)\geq -h\}\}\] {which, after multiplying both sides by $r^{2-2\gamma}$ becomes}
\[\left\{(r,\theta)\,|\,\frac{1}{2}\left(\lambda_2r^{2}\cos^2(\theta)+\lambda_1r^{2} \sin^2(\theta)\right)+
r^{2-\gamma(\nu+2)}+{c}+hr^{2-2\gamma}\geq 0\right\} \]
{The Hill region for the energy levels below, at, and above that of the equilibrium points  $\mathscr{E}_1, \mathscr{E}_2$, is shown in Fig.}~\ref{fig:Hill_region_0}.

\begin{figure}
$\begin{array}{ccc}
\includegraphics[width=0.32\textwidth]{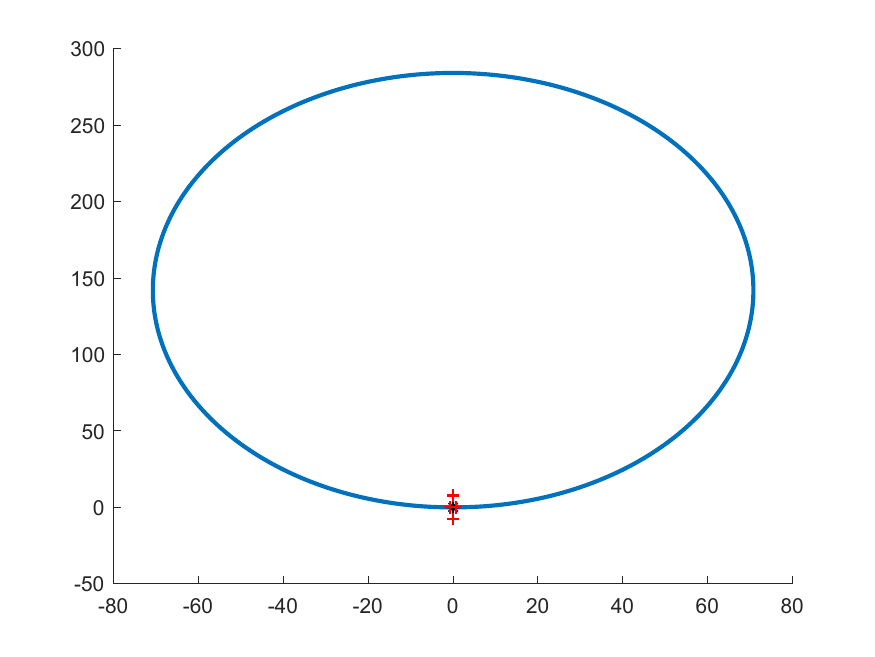}&
\includegraphics[width=0.32\textwidth]{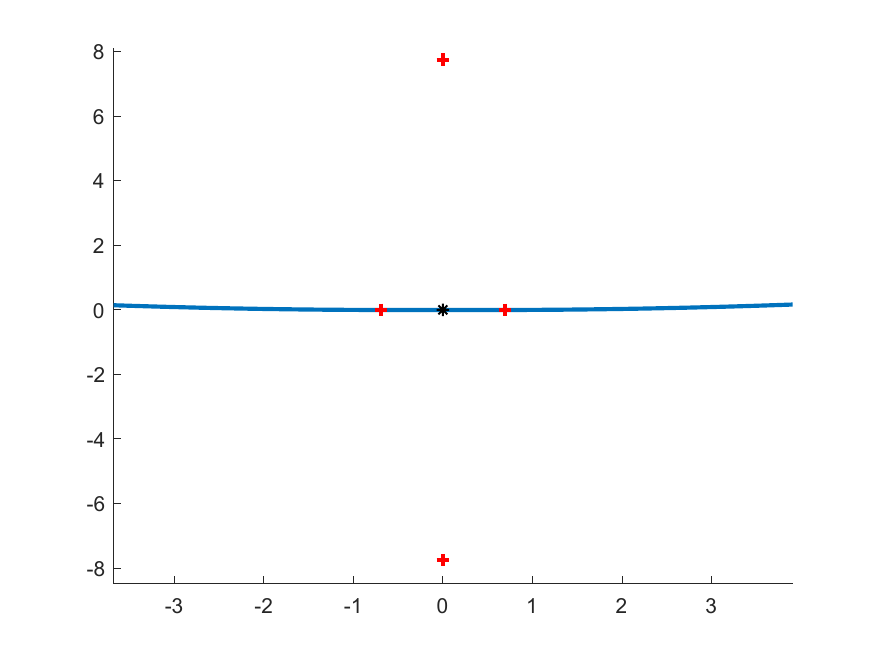}&
\includegraphics[width=0.37\textwidth]{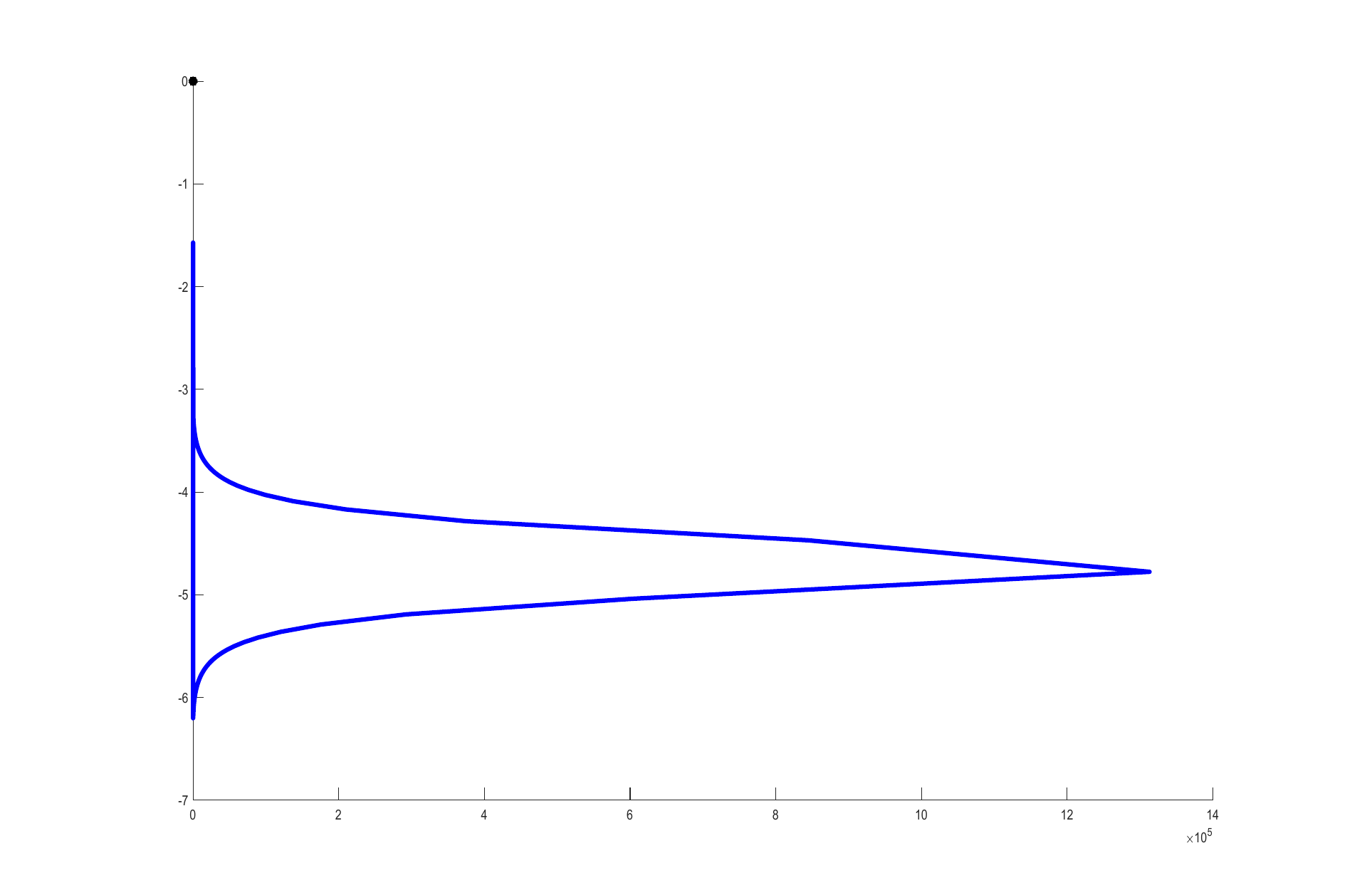}
\\
(a) & (b) &(c)
\end{array}$
\caption{Short period orbit near collision:  (a) in $(x,y)$-Cartesian coordinates, (b) in $(x,y)$-Cartesian coordinates -- magnification near the tertiary (marked by *) and the $\mathscr{E}_1$, $\mathscr{E}_2$, $\mathscr{E}_3$, $\mathscr{E}_4$ equilibrium points (marked by +),  (c) in $(r,\theta)$-coordinates. }
\label{fig:short}\end{figure}

{The system} \eqref{eqn:ham_eqn_2}  {allows the study of the dynamics both near and far from collisions. In particular, it can be used to compute families of orbits that start far from collision and tend asymptotically to collision. For example, we can compute the so called `long' and `short' period families of periodic orbits  near $\mathscr{E}_3, \mathscr{E}_4$, which were studied in} \cite{burgos2016families}. {Such families of periodic orbits were originally considered in}  \cite{deprit1967trojan} {in the context of the planar circular restricted three-body problem, where they emanate from the center-center equilibrium points $L_3$ and $L_4$.
Such equilibrium points do not exist in the Hill three-body problem, but they appear in the Hill four-body problem, as noted in} \cite{burgos2015hill}.
{The long period family of orbits undergoes a bifurcation with the short period family, and the short family approaches a collision with the tertiary as the energy $h$ tends to $+\infty$. An orbit from the short period family, computed in both Cartesian and McGehee coordinates, is shown in Fig.}~\ref{fig:short}.

\section{Collision manifold}
\label{sec:collision_manifold}
From \eqref{eqn:h_r_0}, the intersection between the energy manifold $\mathbf{M}_h$ and the $3$-dimensional hyperplane \[\mathbf{Z}=\{r=0\}\] is a $2$-dimensional manifold corresponding to collisions
\begin{equation}\label{eqn:N_h}
\mathbf{N}=\{H=h\} \cap\{r=0\}.
\end{equation}
It is referred to as   the \emph{collision manifold}. Thus, from \eqref{eqn:h_r_0} we obtain that
\begin{equation}\label{eqn:collision}
\mathbf{N}=\{(r,\theta,v,w)\,|\,r=0,\theta\in\mathbb{T}^1,v^2+w^2=2c\},
\end{equation}
so the collision manifold is a \emph{$2$-dimensional torus} provided $c>0$. Note that this torus is
independent of the energy level $h$, {and is the boundary of each energy manifold $\mathbf{M}_h$.}

The collision manifold is an isolated invariant set for the flow of \eqref{eqn:ham_eqn_2}. If a trajectory
approaches the singularity, i.e., $r\to 0$  as $t\to \pm t^*$, then in the $(r,\theta,v,w)$ coordinates the trajectory approaches the collision manifold $\mathbf{N}$ as $\tau(t)\to \mp \infty$. The argument is the same as in \cite{mcgehee1981double}.

\begin{figure}
\includegraphics[width=0.5\textwidth]{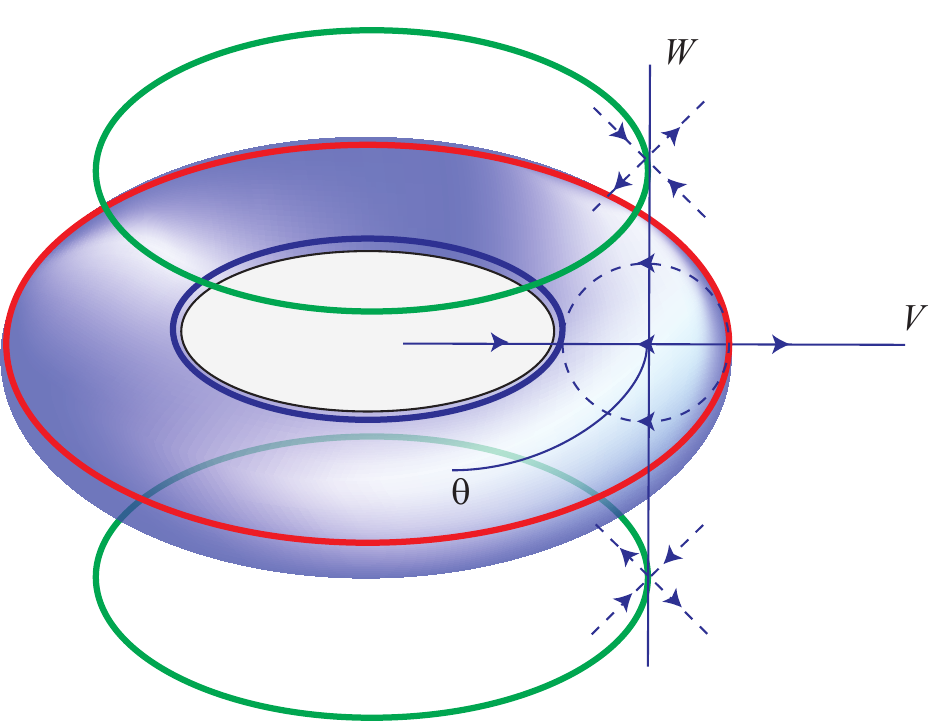}
\caption{The dynamics on $\mathbf{N}$ and $\mathbf{Z}$}\label{fig:collision_manifold}
\end{figure}

Since $r'=0$ when $r=0$, it follows that the set $\mathbf{Z}$ is invariant under the solutions of the system \eqref{eqn:ham_eqn_2}.
Thus, we can consider the restriction of \eqref{eqn:ham_eqn_2} to $\mathbf{Z}$, which is given by
\begin{equation}\label{eqn:ham_eqn_3}\begin{split}
\theta'=&w,\\
v'=&\beta v^2+w^2-\alpha c,\\
w'=&(\beta-1)vw.
\end{split}
\end{equation}
The dynamics on $\mathbf{Z}$ is the skew product between the dynamics in the variables $(v,w)$ and the dynamics in $\theta$.
See Fig. \ref{fig:collision_manifold}.
The solution of the equation in $\theta$ is determined  by the solutions of the $(v,w)$-subsystem, which is independent of $\theta$.
We refer to the  $(v,w)$-subsystem of \eqref{eqn:ham_eqn_3} as the \emph{reduced system associated to the collision manifold}.

Define
\begin{equation}\label{eqn:new_ham}
  K=|w|^\alpha|v^2+w^2-2c|^{1-\beta}.
\end{equation}

We claim that $K$ is an integral of motion for the $(v,w)$-subsystem of \eqref{eqn:ham_eqn_3}. Indeed, using \eqref{eqn:ham_eqn_3} we obtain
\begin{equation*}\begin{split}
  K'=&\alpha|w|^{\alpha-1}w'|v^2+w^2-2c|^{1-\beta}+|w|^\alpha|v^2+w^2-2c|^{-\beta}(2vv'+2ww')\\
  =& \alpha|w|^{\alpha-1}(\beta-1)vw|v^2+w^2-2c|^{1-\beta}\\
  &+|w|^\alpha|v^2+w^2-2c|^{-\beta}(2v(\beta v^2+w^2-\alpha c) +2w(\beta-1)vw)\\
  =&|w|^{\alpha-1}|v^2+w^2-2c|^{-\beta}(\beta-1)vw
  \left[\alpha(v^2+w^2-2c)\right.\\
  &\left.\qquad\qquad\qquad\qquad\qquad\qquad\qquad\quad
   -2(\beta v^2+w^2-\alpha c)-2(\beta-1)w^2\right]\\=&0.
\end{split}\end{equation*}
By \eqref{eqn:collision},  the collision manifold $\mathbf{N}$ intersects the $(v,w)$-plane along the $0$-level set of the integral $K$.

We now describe the geometry of the $(v,w)$-subsystem.
The equilibrium points are  $S_{\pm}=(\pm\sqrt{2c},0)$ and $Q_{\pm}=(0,\pm\sqrt{\alpha c})$.
The circle \[\mathbf{C}=\{v^2+w^2=2c\}\] is invariant under the flow of the subsystem, and passes through the points $S_{\pm}$.
{Thus, $S_{\pm}$ correspond to points on the collision manifold $\mathbf{N}$, while $Q_{\pm}$ do not.}

{The circle $\mathbf{C}$ in the $(v,w)$-plane   corresponds to the collision manifold $\mathbf{N}$,
while the other orbits of the $(v,w)$-subsystem represent projections of orbits on various energy levels onto  the $(v,w)$-plane.}

\begin{figure}
$\begin{array}{cc}
\includegraphics[width=0.5\textwidth]{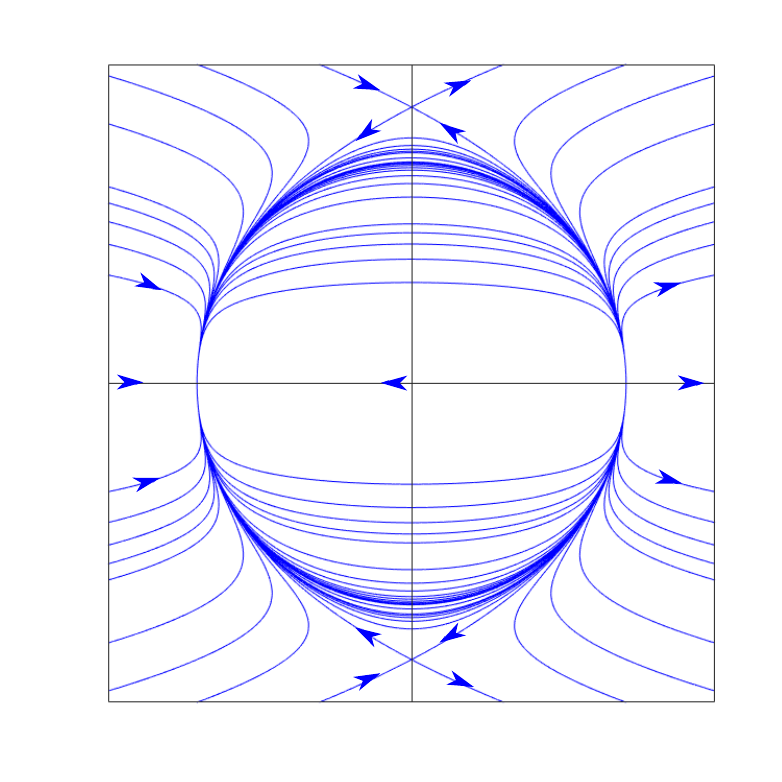}&
\includegraphics[width=0.5\textwidth]{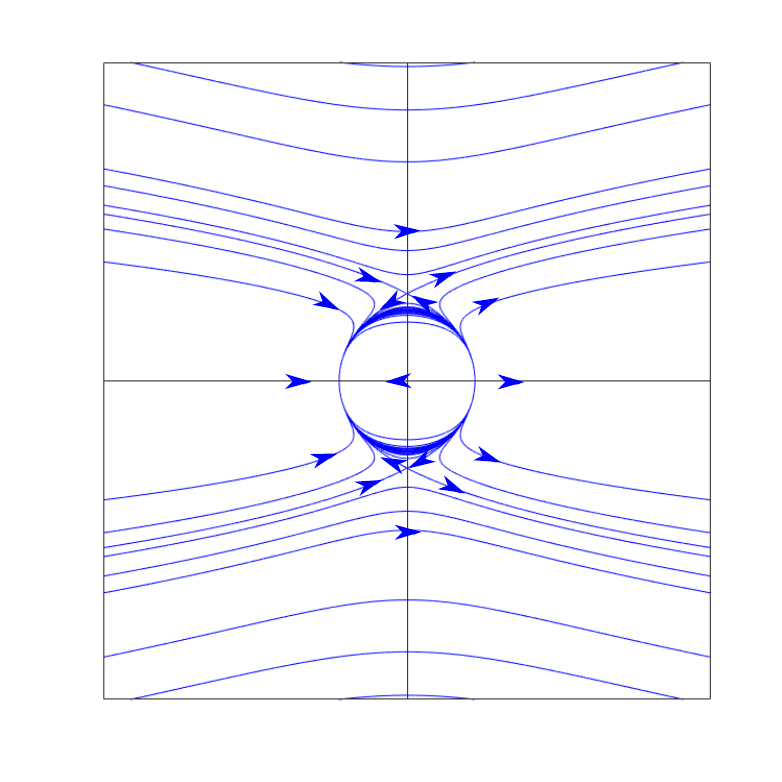}
\end{array}$
\caption{Left: $c=1$. Right: $c=0.1$.}
\label{fig:contour_c_positive}
\end{figure}

The eigenvalues of the linearized system at $Q_{\pm}$ are \[\pm\sqrt{2(\beta-1)\alpha c},\] and since one is positive and the other is negative,  both points are saddle  points.
The eigenvalues of the linearized system at $S_{\pm}$ are
\[\pm 2\beta \sqrt{2c},\quad \pm(\beta-1)\sqrt{2c}.\]
Both eigenvalues at $S_{+}$ are positive hence this is a source.
Both eigenvalues at $S_{-}$ are negative hence this is a sink.
The line $w=0$ is also invariant under the flow, where $v'<0$ for $|v|<\sqrt{2c}$ and $v'>0$ for $|v|>\sqrt{2c}$. The phase portrait is shown in Fig.~\ref{fig:contour_c_positive}.

{Each  point $Q_{\pm}$ has $1$-dimensional stable and unstable manifolds in the  $(v,w)$-plane; these manifolds are asymptotic to $S_\pm$.
In the full phase space, the points $Q_{\pm}$ lie on  circular orbits $\mathscr{C}_{\pm}$ given by}
\[r=0,\quad \theta=\theta_0\pm t \sqrt{\alpha c},\quad v=0, \quad w=\pm\sqrt{\alpha c},\textrm { for } t\in\mathbb{R},\]
{of  energy $+\infty$. %These are limits of direct and retrograde periodic orbits around the tertiary as $h\to +\infty$.
Each circle  $\mathscr{C}_{\pm}$ has $2$-dimensional stable and unstable manifolds in $\mathbf{Z}$.
%??? they do/ do not coincide off the collision manifold.
}

{The points
$S_{\pm}$ lie on the circles of equilibria $\mathscr{E}_{\pm}$ contained in $\mathbf{N}$.
%We stress that only the points $\mathscr{E}_{\pm}$ lie on the collision manifold.
The  circle  $\mathscr{E}_+$  has a $2$-dimensional unstable manifold in $\mathbf{N}$, while the circle $\mathscr{E}_-$ has a $2$-dimensional stable manifold in $\mathbf{N}$;  the stable and unstable manifolds coincide.
In $\mathbf{Z}$  the circle $\mathscr{E}_{+}$ has a $3$-dimensional unstable manifold, and the circle $\mathscr{E}_{-}$ has a $3$-dimensional stable manifold; these manifolds coincide as well.   }

{We summarize the type of orbits that appear near collision:}
\begin{description}
\item[Orbits  beginning and ending in collision] {These orbits form an open set in the phase space, representing the branch of the unstable manifold of $\mathscr{E}_{+}$ that coincides with a branch of stable manifold of $\mathscr{E}_{-}$. Such orbits correspond to initial conditions whose projection onto the $(v,w)$-plane is in $\{(v,w)\,|\, v^2+w^2<2c\}$.}
\item[Orbits that only begin or only end in collision] {These orbits form  open sets in the phase space, representing the branches of the unstable manifold of $\mathscr{E}_{+}$ and of the stable manifold of $\mathscr{E}_{-}$, respectively, whose projection onto the $(v,w)$-plane is in $\{(v,w)\,|\, v^2+w^2>2c\}$.}
\item[Asymptotic orbits  than begin or end in collision] {These orbits represent
branches of the  stable and unstable manifolds of the hyperbolic invariant circles $\mathscr{C}_\pm$.}
\item[Swing-by orbits] {These are orbits coming from afar, passing near the hyperbolic invariant circles $\mathscr{C}_\pm$, and then moving away.}
\end{description}
%\marginpar{Check this}

Recall that for the system \ref{eqn:hill_rotated_planar} we have $\nu=1$, $\alpha=3$, $\beta=\frac{3}{2}$, and $\gamma=\frac{2}{5}$.
By Theorem \ref{thm:branch_2} it follows that each collision solution is branch regularizable.
Since $p=2$ is even,  each extension solution is a `reflection'.

By examining Fig.~\ref{fig:contour_c_positive} we observe that the  collision manifold $\mathbf{N}$ is not an isolated invariant set, and therefore it is not block regularizable. This agrees with the case of the potential \eqref{eqn:potential alpha}, where for $\beta\geq 1$ the collision manifold $\mathbf{N}$ is not an isolated set.

\begin{figure}
$\begin{array}{cc}
\includegraphics[width=0.5\textwidth]{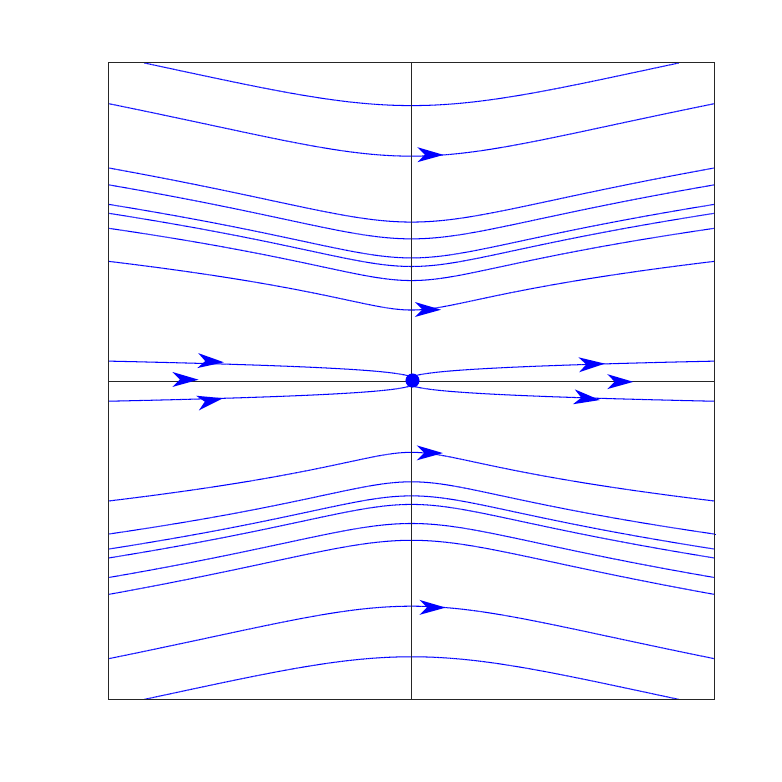}&
\includegraphics[width=0.5\textwidth]{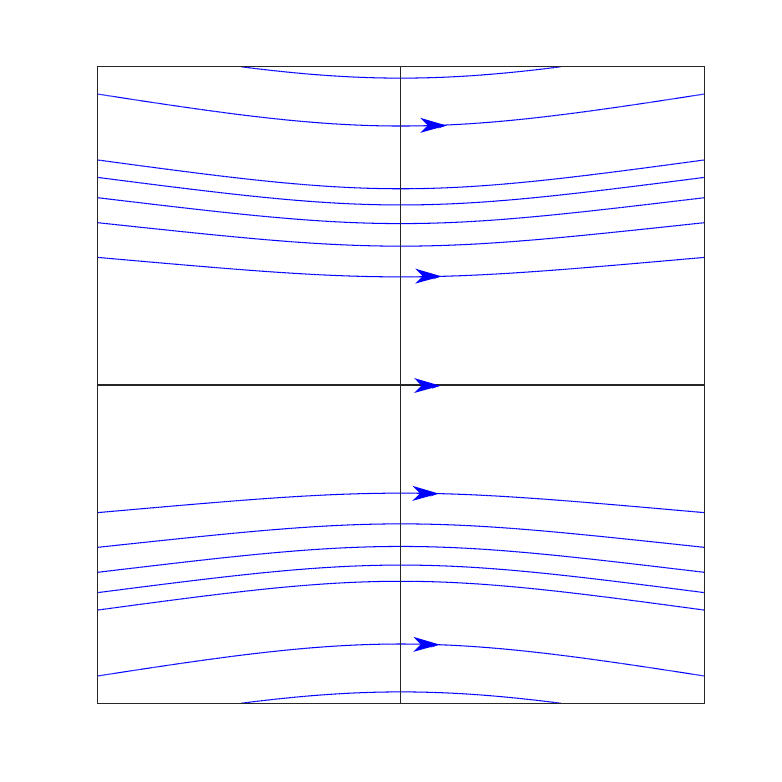}\\
(a) & (b)
\end{array}$
\caption{Phase portrait of the $(v,w)$-subsystem: (a)  $c=0$, (b) $c=-1$.}
\label{fig:c_negative}
\end{figure}

As $c\to 0$, the two saddles, the source, and the sink coalesce through a  double saddle-node bifurcation.
See Fig.~\ref{fig:c_negative}.

For $c=0$ the collision manifold is reduced to a point, and it is both branch and block regularizable.

We now discuss the case when $c<0$.
This describes a situation when the tertiary is a prolate body,
In this case the set of $(v,w)$
with $v^2+w^2=2c$ is the empty set. Thus the collision set $\mathbf{N}$ is empty.
Then the  $(v,w)$-subsystem
\begin{equation}\label{eqn:ham_eqn_4}\begin{split}
v'=&\beta v^2+w^2-\alpha c,\\
w'=&(\beta-1)vw,
\end{split}
\end{equation}
has the property that $v'>0$. The phase portrait is shown in Fig.~\ref{fig:c_negative}.
In this case there are no collisions.

The physical interpretation is the following. Denoting $c=-\tilde{c}$ where $\tilde{c}>0$, the Hamiltonian \eqref{eqn:hill_rotated_planar} becomes
\begin{equation}\label{eqn:hill_repel}\begin{split}
H=&\frac{1}{2}(y_1^2+y_2^2)+x_2 y_1-x_1y_2+Ax_1^2+B x_2^2\\
&\quad  -\frac{1}{|x|^{\nu/2}}+ \frac{\tilde{c}}{|x|^{\alpha/2}}.
\end{split}\end{equation}
The term $-\frac{1}{|x|^{\nu/2}}$ in the potential corresponds to an attractive force, and the term $\frac{\tilde{c}}{|x|^{\alpha/2}}$ corresponds to a repulsive force. When the particle approaches the tertiary, since $1\leq \nu<\alpha$ the repulsive force becomes dominating, preventing collisions between the particle and the tertiary to occur.
This situation is also described in \cite{saari1974regularization}.

\begin{rem}
{One can  consider a simple model that takes into account the size of the asteroid}. {Since $\nu=1$, $\alpha=3$, $\beta=\frac{3}{2}$, and $\gamma=\frac{2}{5}$, the powers of $r$ that appear in } \eqref{eqn:ham_eqn_2}  {are}
\[r^{2-\gamma(1+\nu)}=r^{4/5}\gg r\gg r^2.\]
{We can neglect the  powers $r^k$ of $r$ with $k>4/5$. Then collisions correspond to setting $r=R_3$, the average radius of the tertiary;
in the case of Hektor, in the normalized units $R_3=1.18716\times 10^{-7}$.
Then } \eqref{eqn:ham_eqn_2} {yields}
\begin{equation}\label{eqn:ham_eqn_5}\begin{split}
\theta'=&w,\\
v'=&( \beta v^2+w^2-\alpha c)-\nu R_3^{2-\gamma(\nu+2)},\\
w'=&(\beta-1)vw.
\end{split}
\end{equation}
{This system is essentially the same as the system} \eqref{eqn:ham_eqn_3} {with the term $-\alpha c$ replaced with the term $-\alpha c-\nu R_3^{2-\gamma(\nu+2)}$.
Then the analysis of collisions is similar to the one above.
Another possibility could be to neglect the  powers $r^k$ of $r$ with $k>1$.
Of course, it is possible to consider more sophisticated models that take into account the dumbbell  shape of Hektor or more general asteroid shapes,  in which case the  gravitational potential} \eqref{eqn:potential} {needs to be modeled differently, e.g.,} \cite{lam2021surface}.
\end{rem}

\begin{rem}
There are several  moons in the Solar System that are considered to be approximately prolate spheroids in shape, for example, Uranus' moons Cordelia, Cressida, Desdemona, Juliet, Ophelia, and Rosalind.
\end{rem}

\begin{rem}
 {The  $(v,w)$-subsystem of} \eqref{eqn:ham_eqn_3} {is undergoing another  bifurcation at $\alpha=2$ when  $c>0$ is held fixed.
When $\alpha=2$ the points $Q_{\pm}$ lie on the collision manifold.
When  $\alpha<2$, the  points  $Q_{\pm}$ become centers, and the points $S_\pm$  become saddles. The phase portrait is as in Fig.}~\ref{fig:bifurcation_alpha}.
\begin{figure}
\includegraphics[width=0.5\textwidth]{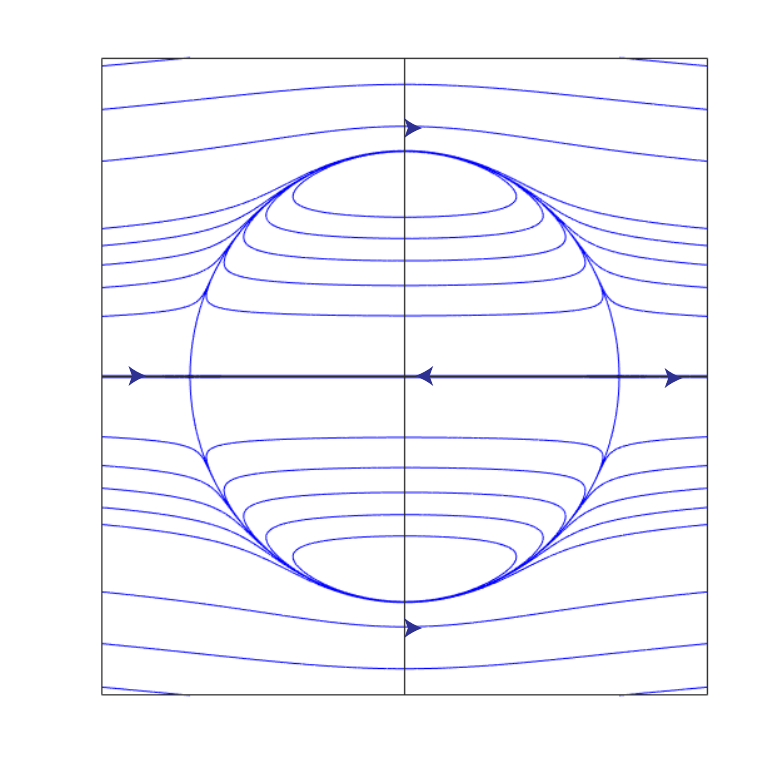}
\caption{Phase portrait of the $(v,w)$-subsystem for $\alpha<2$ and $c>0$.}
\label{fig:bifurcation_alpha}\end{figure}
\end{rem}

\begin{rem}
In the case when $c=0$, the term $-\frac{c}{|x|^\alpha}$  in \eqref{eqn:hill_rotated_simplified} vanishes.
Then one can perform the coordinate change \eqref{eqn:mcgehee} with $\gamma=\frac{2}{2+\nu}$ and $\beta=\frac{\nu}{2}$, as in \cite{mcgehee1981double}.
{The resulting collision manifold is a torus which intersects the $(v,w)$-plane in a circle as in} Fig.~\ref{fig:McGehee}. {Note that the phase portrait is qualitatively the same as in  Fig.}~\ref{fig:bifurcation_alpha}.
{It contains  two circles of equilibria located at $v =\pm\sqrt{2}$ and a cylinder of
orbits given by $w=0$ connecting the two circles.
%In the Hill three-body problem these two circles
}
The collision set  is both branch and block regularizable.
It is interesting that this coordinate change leads to a different collision manifold from the one in \eqref{eqn:collision}, but nevertheless its branch and block regularization properties are the same.
\begin{figure}
\includegraphics[width=0.5\textwidth]{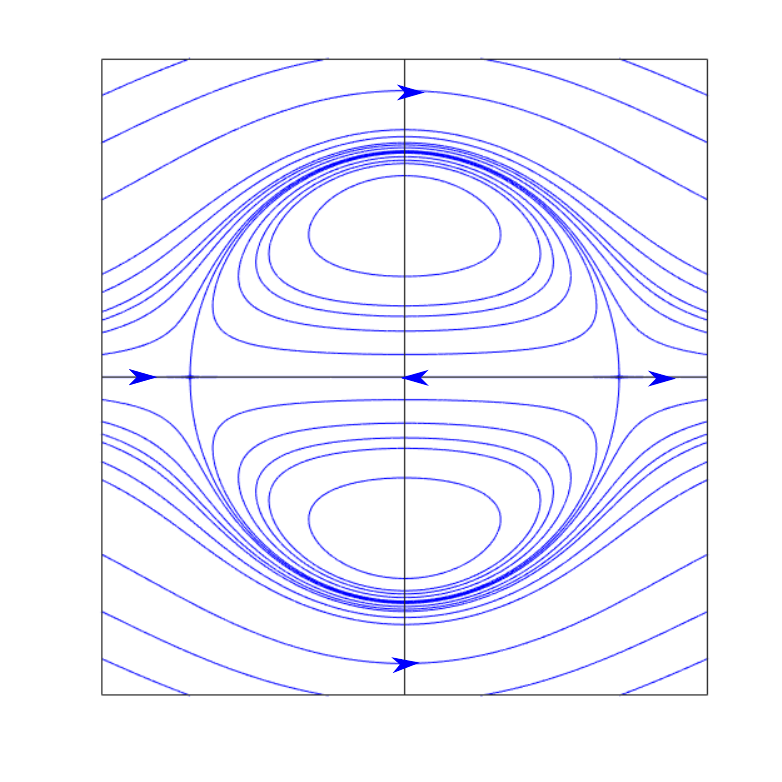}
\caption{Phase portrait of the $(v,w)$-subsystem  for $c=0$ with the coordinate change \eqref{eqn:mcgehee} with $\gamma=\frac{2}{2+\nu}$ and $\beta=\frac{\nu}{2}$.}
\label{fig:McGehee}\end{figure}
\end{rem}

\section{Conclusions}
In this paper we provide an explicit McGehee coordinate transformation to regularize collision in the planar Hill four-body problem with oblate bodies.
This transformation can be used to understand the behavior of collision and near-collision orbits. In particular, our formulas can be implemented in numerical integrators to compute  orbits that pass close to an oblate Jupiter's trojan asteroid.

We also describe the collision manifold and show that it undergoes a bifurcation as the oblateness coefficient of the asteroid passes through the zero value.
{We note here that the bifurcation observed for this system is very different from the one described by} \cite{mcgehee1981double} {for the  potential energy} $U({\bf x})=|{\bf x}|^{-\alpha}$ {in}  \eqref{eqn:potential alpha}, {which undergoes a bifurcation when  the parameter $\alpha$ passes through the critical value} $\alpha_{\textrm{cr}}=2$.

It is interesting to note that when the oblateness approaches  zero (and hence the gravitational potential becomes Newtonian), the limiting collision manifold that we obtain is not the same as the collision manifold obtained by applying the McGehee coordinate transformation to the Newtonian potential. It would be interesting to see if there is a  McGehee-type  coordinate transformation for which the limiting collision manifold is the same as in the Newtonian case.
Another interesting problem would be to extend these results to the spatial Hill four-body problem with oblate bodies.

\section*{Acknowledgements} We are grateful to Jaime Burgos-Garc\'ia for useful discussions.
\bibliographystyle{alpha}
\bibliography{regularization}

\end{document}